# Binomial expansion of Newton's method

Shunji Horiguchi

**Introduction.**

In 1673, Yoshimasu Murase made a cubic equation to obtain the thickness of a hearth. He introduced two kinds of recurrence formulas of square $x_k^2$ and the deformation. We find that the three formulas lead to a Horner's method and extension of a Newton's method at the same time(Ref.[1],[2]). This shows originality of Wasan(Japanese native mathematics developed in Japan in the Edo era(1603-1868:national isolation)). Furthermore, in this paper, we show that the extended Newton's method leads to the binomial expansion of Newton's method that the convergences become the quadratic and linearly. Next we give convergence comparisons of the binomial expansion of Newton's method and Newton's method using the curvature and convex-concave of curve. Finally, we give examples of the numerical calculations.

We start with a change of variable of the function.

**1. Function $y=g(t)$ defined by $x=t^{1/q}$ of $y=f(x)$.**

**Definition 1.1.** Let $x=t^{1/q}$ where $q$ is a real number that is not 0. We define the function $g(t)$ such as

$$g(t) := f(t^{1/q}) = f(x). \tag{1.1}$$

Because $g(x^q)=f(x)$, the graph of $g(x)$ is extended and contracted by $x^q=t$ in the $x$-axis, without changing the height of $f(x)$. Expansion and contraction come to object in $|x|<1$ and $|x|>1$.

**Lemma 1.2.** $g'(x^q)$, $g''(x^q)$ are represented by $f'(x)$, $f''(x)$ as follows.

$$g'(x^q) = \frac{f'(x)}{qx^{q-1}} \tag{1.2}$$

$$g''(x^q) = \begin{cases} \dfrac{xf''(x)+(1-q)f'(x)}{q^2 x^{2q-1}} & \tag{1.3} \\[2ex] \dfrac{f''(x)\left(1+\dfrac{1-q}{x}\dfrac{f'(x)}{f''(x)}\right)}{\left(qx^{q-1}\right)^2}, \; f''(x) \neq 0 & \tag{1.4} \end{cases}$$

From this lemma we get the next theorem.

**Theorem 1.3.** The curvature of the curve $y=g(x)$ at the point $x^q$ is this.

$$\mu_q(t) = \frac{g''(t)}{(1+g'(t)^2)^{3/2}} = \mu_q(x^q) = \begin{cases} \dfrac{f''(x)+\dfrac{(1-q)f'(x)}{x}}{\left(qx^{q-1}\right)^2\left(1+\left(\dfrac{f'(x)}{qx^{q-1}}\right)^2\right)^{3/2}} & \tag{1.5} \\[3ex] \dfrac{f''(x)\left(1+\dfrac{1-q}{x}\dfrac{f'(x)}{f''(x)}\right)}{\left(qx^{q-1}\right)^2\left(1+\left(\dfrac{f'(x)}{qx^{q-1}}\right)^2\right)^{3/2}}, \; f''(x) \neq 0 & \tag{1.6} \end{cases}$$

These become $\mu(x) = f''(x)/(1+f'(x)^2)^{3/2}$ of $f(x)$ if $q=1$ in particular.

## 2. Extension of Newton's method and binomial expansion of Newton's method.

**Definition 2.1.** The recurrence formula to approximate a root of the equation $f(x)=0$

$$x_{k+1} = x_{k+1} - \frac{f(x_k)}{f'(x_k)}, \quad (k=0,1,2,\cdots) \tag{2.1}$$

is called Newton's method($N$-method,1669) or Newton-Raphson's method(1690).

Newton's method is a method of giving the initial value $x_0$, calculating $x_1, x_2, \cdots$ one after another, and to determine for a root.

Applying the Newton's method to $g(t)$, we have

$$t_{k+1} = t_k - \frac{g(t_k)}{g'(t_k)}, \quad t_{k+1} = t_k - \frac{f(t_k^{1/q})}{f'(t_k^{1/q}) \tfrac{1}{q} t_k^{1/q-1}}. \tag{2.2}$$

This means the intersection $t_{k+1} = x_{k+1}^q$ with the $t(x)$-axis of the tangent in the point $(t_k, g(t_k)) = (x_k^q, g(x_k^q))$ of the graph of $y=g(t)(g(x))$. Returning to the variable $x$ by $x^q = t$, we get an extension of Newton's method below.

**Definition 2.2.** For equation $f(x)=0$, we call the next recurrence formulas the extended Newton's method($EN$-method) or Tsuchikura-Horiguchi's method($TH$-method).

$$x_{k+1}^q = x_k^q - qx_k^{q-1}\frac{f(x_k)}{f'(x_k)} \; (q \neq 0, q \in R) \iff x_{k+1} = \left[ x_k^q - qx_k^{q-1}\frac{f(x_k)}{f'(x_k)} \right]^{\frac{1}{q}} \tag{2.3}$$

Here, if $q=1$ then the formula (2.3) becomes Newton's method.

**Proposition 2.3.** If $\alpha$ is a simple root($m(>1)$ multiple root resp.) of $f(x)=0$, then $\alpha^q$ becomes the simple root($m$ multiple root resp.) of $g(x)$.

When $q$ is an integer greater than or equal to 2, the extended $N$-method (2.3) is the first term +second term of the $q$th power of $N$-method (2.1). Therefore, we extend the $N$-method by the binomial expansion. First, we give Newton's general binomial coefficient in 1665.

**Definition 2.4.** The following formula is called Newton's general binomial coefficient.

$$\binom{r}{i} = \frac{r(r-1)(r-2)\cdots(r-i+1)}{i!}, \; r : \text{real number} \tag{2.4}$$

**Definition 2.5.** Let $q(\neq 0)$ be a real number. The following formula is called the binomial expansion from the first term to the $m+1$ term of the $q$th power of Newton's method.

$$x_{k+1}^q = \left[ x_k - \frac{f(x_k)}{f'(x_k)} \right]^{q,m} := \sum_{i=0}^{m} \binom{q}{i} x_k^{q-i} \left[ -\frac{f(x_k)}{f'(x_k)} \right]^i$$

$$= x_k^q - qx_k^{q-1}\frac{f(x_k)}{f'(x_k)} + \frac{q(q-1)}{2!} x_k^{q-2}\left[ -\frac{f(x_k)}{f'(x_k)} \right]^2 + \cdots + \frac{q(q-1)\cdots(q-m+1)}{m!} x_k^{q-m}\left[ -\frac{f(x_k)}{f'(x_k)} \right]^m \tag{2.5}$$

Especially if $m=1$ ( $m \geq q$(integer) resp.) then formula (2.5) becomes $EN$-method($N$-method resp.).

## 3. The convergences of binomial expansion of Newton's method.

**Lemma 3.1.** In the sequence $\{x_n\}$, let $\lim_{n\to\infty} x_n = \alpha$ and $q, r$ an arbitrary real constant that is not 0, respectively. In this case, following formula holds for large enough integer $n$.

$$x_n^q - \alpha^q \doteqdot \frac{q}{r}\alpha^{q-r}(x_n^r - \alpha^r) \tag{3.1}$$

Proof. Applying L'Hospital's rule to $(x^q - \alpha^q)/(x^r - \alpha^r)$, (3.1) is obtained. □

From now, it is assumed that the initial value $x_0$ is close to the root $\alpha$ of $f(x)=0$.

**Theorem 3.2.** Let $\alpha$ be a simple root for $f(x)=0$ i.e., $f'(\alpha) \neq 0$. Then Newton's method to the quadratic convergence of the following formula.

$$x_{k+1} - \alpha \doteqdot \frac{1}{2}\frac{f''(\alpha)}{f'(\alpha)}(x_k - \alpha)^2 \tag{3.2}$$

If $\alpha$ is $m(\geq 2)$ multiple root, then it will become the linearly convergence of the following formula.

$$x_{k+1} - \alpha \doteqdot \left(1 - \frac{1}{m}\right)(x_k - \alpha) \tag{3.3}$$

Proof should see at the books of standard numerical computation.

**Remark.** Concerning choosing the initial value $x_0$, the number of iterations until $x_k$ converges on a root change. Moreover, it may not be converged on a root.

**Theorem 3.3.** Let $\alpha(\neq 0)$ be a simple root for $f(x)=0$ i.e., $f'(\alpha) \neq 0$. The binomial expansion (2.5) from the first term to the second term or more is the next quadratic convergence.

$$x_{k+1} - \alpha \doteqdot \frac{1}{2}\left[\frac{f''(\alpha)}{f'(\alpha)} + \frac{1-q}{\alpha}\right](x_k - \alpha)^2 \tag{3.4}$$

In particular, if $m \geq q$ (integer) in (2.5), then the convergence of (2.5) becomes (3.2).

If $\alpha$ is $m(\geq 2)$ multiple root, then it will become linearly convergence of the following formula.

$$x_{k+1} - \alpha \doteqdot \left(1 - \frac{1}{m}\right)(x_k - \alpha) \tag{3.5}$$

Proof. In case of $\alpha$ is a simple root.

(i) Convergence from the first term to the second term of (2.5).

In this case, convergence of (2.5) is that *EN*-method (2.3). Because $\alpha^q$ is a simple root for $g(t)=0$, Newton's method for $g(t)$ becomes the quadratic convergence of the following formula.

$$t_{k+1} - \alpha^q \doteqdot \frac{1}{2}\frac{g''(\alpha^q)}{g'(\alpha^q)}(t_k - \alpha^q)^2 \tag{3.6}$$

Here by the (1.2) and (1.3), substituting $g'(\alpha^q), g''(\alpha^q)$ into (3.6) gives the next formula.

$$x_{k+1}^q - \alpha^q \doteqdot \frac{1}{2}\frac{\alpha f''(\alpha) + (1-q)f'(\alpha)}{q^2 \alpha^{2q-1}}\frac{q\alpha^{q-1}}{f'(\alpha)}(x_k^q - \alpha^q)^2. \tag{3.7}$$

Here by the (3.1), next formula is obtained.

$$q\alpha^{q-1}(x_{k+1}-\alpha) \doteqdot \frac{1}{2}\frac{\alpha f''(\alpha)+(1-q)f'(\alpha)}{q\alpha^q}\frac{1}{f'(\alpha)}q^2\alpha^{2q-2}(x_k-\alpha)^2 \qquad (3.8)$$

(ii) Convergence from the first term to the third term or more of (2.5).
Substituting (3.9) into (2.5) and omitting the terms of $x_k-\alpha$ cubed or higher give the following formula (3.10).

$$\frac{f(x_k)}{f'(x_k)} \doteqdot x_k - \alpha - \frac{1}{2}\frac{f''(\alpha)}{f'(\alpha)}(x_k-\alpha)^2 \qquad (3.9)$$

$$x_{k+1}^q - \alpha^q \doteqdot x_k^q - \alpha^q - qx_k^{q-1}\left[x_k-\alpha-\frac{1}{2}\frac{f''(\alpha)}{f'(\alpha)}(x_k-\alpha)^2\right] + \frac{q(q-1)}{2!}x_k^{q-2}\left[x_k-\alpha-\frac{1}{2}\frac{f''(\alpha)}{f'(\alpha)}(x_k-\alpha)^2\right]^2 \qquad (3.10)$$

Here by (3.1), we obtain the following formula from (3.11).

$$x_{k+1}-\alpha \doteqdot \left[-\frac{1}{\alpha^{q-1}}(q-1)\alpha^{q-2}+\frac{1}{2}\frac{f''(\alpha)}{f'(\alpha)}+\frac{1}{\alpha}\frac{q-1}{2}\right](x_k-\alpha)^2 = \left[\frac{1}{2}\frac{f''(\alpha)}{f'(\alpha)}+\frac{1}{2}\frac{1-q}{\alpha}\right](x_k-\alpha)^2 \qquad (3.11)$$

In case of $\alpha$ is $m(\geq 2)$ multiple roots.
In this case $f(x)$ is expressed by the following formula.

$$f(x) = (x-\alpha)^m g(x), \quad g(\alpha) \neq 0 \qquad (3.12)$$

$$\frac{f(x)}{f'(x)} = \frac{(x-\alpha)^m g(x)}{m(x-\alpha)^{m-1}g(x)+(x-\alpha)^m g'(x)} = \frac{(x-\alpha)g(x)}{mg(x)+(x-\alpha)g'(x)}$$

$$= \frac{1}{m}(x-\alpha) - \frac{g'(x)}{m}\frac{1}{mg(x)+(x-\alpha)g'(x)}(x-\alpha)^2 \qquad (3.13)$$

Substituting above formula into (2.5) and omitting the terms of the square or more of $x_k-\alpha$ give (3.5). □

**4. Varieties of formulas to compare the convergences for the binomial expansion of Newton's method.**

**Theorem 4.1.** Let $\alpha(\neq 0)$ be a simple root of $f(x)=0$ and $f''(\alpha)\neq 0$. If an appropriate initial value $x_0$ is selected for $q$ that satisfies formula (4.1) then the convergence to $\alpha$ of the $q$th power of binomial expansion of Newton's method is equal to or faster than that Newton's method.

$$\left|1+\frac{f'(\alpha)}{f''(\alpha)}\frac{1-q}{\alpha}\right| \leq 1 \Leftrightarrow 0 \leq \frac{f'(\alpha)}{f''(\alpha)}\frac{q-1}{\alpha} \leq 2 \qquad (4.1)$$

Proof. If we compare the coefficient of $(x_k-\alpha)^2$ in (3.4) and that (3.2) then we get

$$\frac{1}{2}\left|\frac{f''(\alpha)}{f'(\alpha)}+\frac{1-q}{\alpha}\right| \leq \frac{1}{2}\left|\frac{f''(\alpha)}{f'(\alpha)}\right|. \qquad (4.2)$$

The formula (4.1) is obtained from (4.2). □

**Theorem 4.2.** Let $\alpha(\neq 0)$ be a simple root of $f(x)=0$, and $f''(\alpha)=0$ (i.e., the graph of $f(x)$ is nearly the straight line in the neighborhood of the point $\alpha$.). In this case (4.3) holds.

$$|\mu(\alpha)|=0 \leq |\mu_q(\alpha^q)| \, (q \neq 1) \tag{4.3}$$

This is equivalent to the convergence to $\alpha$ of Newton's method equals to or faster than that the $q$th power of binomial expansion of Newton's method.

Proof. By deforming the formula to $\mu_q(\alpha^q)$, we compare it with $\mu(\alpha)$.

$$|\mu(\alpha)|=0 \leq |\mu_q(\alpha^q)|= \frac{\left|0+\dfrac{(1-q)f'(\alpha)}{\alpha}\right|}{\left(q\alpha^{q-1}\right)^2\left(1+\left(\dfrac{f'(\alpha)}{q\alpha^{q-1}}\right)^2\right)^{3/2}} = \frac{\left|\dfrac{(1-q)f'(\alpha)}{\alpha}\right|}{\left(q\alpha^{q-1}\right)^2\left(1+\left(\dfrac{f'(\alpha)}{q\alpha^{q-1}}\right)^2\right)^{3/2}} \Leftrightarrow 0 \leq |1-q|$$

$$0 \leq |1-q| \Leftrightarrow \frac{1}{2}\left|\frac{f''(\alpha)(=0)}{f'(\alpha)}\right| \leq \frac{1}{2}\left|\frac{f''(\alpha)(=0)}{f'(\alpha)}+\frac{1-q}{\alpha}\right| \tag{4.4}$$

We get the conclusion by (4.4). □

Following are the results related to the convex-concave of curve and the formulas for comparing convergences of the binomial expansions of Newton's method.

**Lemma 4.3.** Let $x \neq 0$ and $f''(x) \neq 0$. Then a necessary and sufficient condition for $g''(x^q)$ and $f''(x)$ are the same sign (opposite sign resp.) is

$$1+\frac{f'(x)}{f''(x)}\frac{1-q}{x} > 0 (<0 \text{ resp.}). \tag{4.5}$$

Proof. Because

$$g''(x^q) = \frac{xf''(x)+(1-q)f'(x)}{q^2 x^{2q-1}} = \frac{f''(x)\left(1+\dfrac{(1-q)f'(x)}{xf''(x)}\right)}{(qx^{q-1})^2}, \tag{4.6}$$

according to $1+(1-q)f'(x)/xf''(x) > 0 (<0 \text{ resp.})$, $g''(x^q)$ and $f''(x)$ become the same sign(opposite sign resp.). □

We get the next theorem from lemma 4.3, directly.

**Theorem 4.4.** Let $\alpha(\neq 0)$ be a simple root of $f(x)=0$, and $f''(\alpha) \neq 0$. We divide the formula (4.1) of theorem 4.1 into positive and negative range as follows.

$$-1 \leq 1+\frac{f'(\alpha)}{f''(\alpha)}\frac{1-q}{\alpha} < 0 \tag{4.7}$$

$$0 < 1+\frac{f'(\alpha)}{f''(\alpha)}\frac{1-q}{\alpha} \leq 1 \tag{4.8}$$

If $q$ satisfies the condition (4.8)((4.7) resp.), then the convex-concave of curve of $g(x)$ in the neighborhood of $g(\alpha^q)(=0)$ and that $f(x)$ in the neighborhood of $f(\alpha)(=0)$ are the same(opposite resp.).

**Theorem 4.5.** Let the conditions be the same as the above theorem. If $q$ satisfies the inequality (4.9) then the convergence to $\alpha$ of $q$th power of binomial expansion of Newton's method is equal to or faster than that Newton's method.

$$-\frac{|f''(\alpha)|}{(q\alpha^{q-1})^2} \leq g''(\alpha^q) \leq \frac{|f''(\alpha)|}{(q\alpha^{q-1})^2} \tag{4.9}$$

Proof. By the formula

$$\left|\frac{xf''(x)+(1-q)f'(x)}{q^2 x^{2q-1}}\right|\left|\frac{q^2 x^{2q-1}}{xf''(x)}\right| = \left|1+\frac{f'(x)}{f''(x)}\frac{1-q}{x}\right| \leq 1 \tag{4.10}$$

and (4.1) of theorem 4.1, (4.9) is obtained. □

**Corollary 4.6.** If $|q\alpha^{q-1}|=1$ then inequality (4.9) becomes

$$-|f''(\alpha)| \leq g''(\alpha^q) \leq |f''(\alpha)|. \tag{4.11}$$

The following are the results related to the curvature and the formula (4.1) for comparing the convergences of the binomial expansion of Newton's method.

**Theorem 4.7.** Let $\alpha(\neq 0)$ be a simple root of $f(x)=0$, and $f''(\alpha) \neq 0$. Suppose that the curvature $\mu_q(x)$ of $g(x)$ satisfies the condition

$$|\mu_q(\alpha^q)| \leq \frac{|f''(\alpha)|}{(q\alpha^{q-1})^2 \left(1+\left(\frac{f'(\alpha)}{q\alpha^{q-1}}\right)^2\right)^{3/2}}. \tag{4.12}$$

Then the convergence to $\alpha$ of $q$th power of binomial expansion of Newton's method is equal to or faster than that Newton's method.

Proof. By the formula

$$\frac{f''(\alpha)+\frac{(1-q)f'(\alpha)}{\alpha}}{(q\alpha^{q-1})^2\left(1+\left(\frac{f'(\alpha)}{q\alpha^{q-1}}\right)^2\right)^{3/2}}(q\alpha^{q-1})^2\left(1+\left(\frac{f'(\alpha)}{q\alpha^{q-1}}\right)^2\right)^{3/2}\frac{1}{f''(\alpha)} = 1+\frac{f'(\alpha)}{f''(\alpha)}\frac{1-q}{\alpha} \tag{4.13}$$

and (4.1) of theorem 4.1, (4.12) is obtained. □

**Theorem 4.8.** Let the conditions be same as the above theorem. If

$$\left|\mu_q(\alpha^q)\right| \leq \left|\mu(\alpha)\right| \tag{4.14}$$

and

$$\left(q\alpha^{q-1}\right)^2\left(1+\left(\frac{f'(\alpha)}{q\alpha^{q-1}}\right)^2\right)^{3/2} \leq \left(1+f'(\alpha)^2\right)^{3/2} \tag{4.15}$$

hold, then the convergence to $\alpha$ of the $q$th power of binomial expansion of Newton's method is equal to or faster than that Newton's method.

Proof. From (4.14) we obtain

$$|\mu_q(\alpha^q)| = \frac{\left|f''(\alpha) + \frac{(1-q)f'(\alpha)}{\alpha}\right|}{\left(q\alpha^{q-1}\right)^2\left(1+\left(\frac{f'(\alpha)}{q\alpha^{q-1}}\right)^2\right)^{3/2}} \leq \frac{|f''(\alpha)|}{(1+f'(\alpha)^2)^{3/2}} = |\mu(\alpha)| \qquad (4.16)$$

$$\Leftrightarrow \left|1 + \frac{f'(\alpha)}{f''(\alpha)}\frac{1-q}{\alpha}\right| \leq \frac{\left(q\alpha^{q-1}\right)^2\left(1+\left(\frac{f'(\alpha)}{q\alpha^{q-1}}\right)^2\right)^{3/2}}{(1+f'(\alpha)^2)^{3/2}} \qquad (4.17)$$

Assertion is obtained from (4.17) and (4.1). □

**5. Examples of convergence comparisons of the numerical calculations of Newton's method and binomial expansion of Newton's method.**

Numerical calculations are performed using the standard 10 digits of Microsoft excel.

**Example 5.1.** A quadratic equation

$$f(x) = (x-1)(x-2) = x^2 - 3x + 2 = 0. \qquad (5.1.1)$$

The roots of (5.1.1) are $\alpha = 1, 2$. Because $f'(x) = 2x - 3$, $f''(x) = 2$, condition (4.2) becomes

$$0 \leq \frac{2\alpha - 3}{2}\frac{q-1}{\alpha} \leq 2. \qquad (5.1.2)$$

In case of $\alpha = 1$, (5.1.2) becomes $-3 \leq q \leq 1$.

We examine the speed of convergence of the $q$th power of binomial expansion of Newton's method in $-4 \leq q \leq 2$. The results of the calculations are table 5.1.1. Here the first column represents the real number $q$ and $m$ of formula (2.5). The second, third and fourth column represents the initial value $x_0$, the number of iterations that the $q$th power of binomial expansion of Newton's method to converge to the root 1, the absolute error

|the value 1 of the convergence of the numerical calculation $x_k$ − root 1|

, respectively. We evaluate by the absolute errors when the two convergence times are equal. Thus, these numerical calculations are compatible with the theory of theorem 4.1.

In case of $\alpha = 2$, (5.1.2) becomes $1 \leq q \leq 9$.

The results of the calculations are table 5.1.2. We can also confirm that the numerical calculations are compatible with the theory of theorem 4.1.

**Example 5.2.** Numerical calculations of formulas (4.12),(4.14),(4.15),(4.17) for $f(x) = x^2 - 3x + 2$

In case of $\alpha = 1$, (4.12) of theorem 4.7 becomes

$$|\mu_q(1^q)| = \frac{|q+1|}{q^2\left(1+\frac{1}{q^2}\right)^{3/2}} \leq \frac{2}{q^2\left(1+\frac{1}{q^2}\right)^{3/2}}. \qquad (5.2.1)$$

Indeed, by calculating the left and right sides of (5.2.1) for $q$ in the table 5.2.1, we get the numbers there. For each $q$ in $-3 \leq q \leq 1$, the numbers of the second column and third column satisfy the condition (5.2.1). Formulas (4.14),(4.15) of theorem 4.8 hold for $q = -0.5, 0.5, 1$.

In case of $\alpha=2$, (4.12) of theorem 4.7 becomes

$$|\mu_q(2^q)|=\frac{\left|2+\frac{1-q}{2}\right|}{\left(q2^{q-1}\right)^2\left(1+\left(\frac{1}{q2^{q-1}}\right)^2\right)^{3/2}}\leq\frac{2}{\left(q2^{q-1}\right)^2\left(1+\left(\frac{1}{q2^{q-1}}\right)^2\right)^{3/2}}. \quad (5.2.2)$$

Indeed, by calculating the left and right sides of (5.2.2) for $q$ in table 5.2.2, formula (5.2.2) holds in $1\leq q\leq 9$. Formulas (4.14),(4.15) of theorem 4.8 hold the equal sign only $q=1$.

| Table 5.1.1 Caluculations of (2.5) in case of root 1 | | | |
|---|---|---|---|
| $q$ | $x_0$ | iterations | absolute error |
| N-method | 0.85 | 4 | 6.66134E-15 |
| (q=1) | 1.3 | 5 | 1.67777E-12 |
| -4(2 term) | | 5 | |
| -4(3 term) | 0.85 | 5 | |
| -4(4 term) | | 4 | 2.14051E-12 |
| -3(2 term) | | | 2.29594E-13 |
| -3(3 term) | 1.3 | 5 | 1.06581E-14 |
| -3(4 term) | | | 1.06581E-14 |
| -2(2 term) | | 4 | |
| -2(3 term) | 1.3 | 5 | 7.06102E-14 |
| -2(4 term) | | 5 | 3.26406E-14 |
| -1(2 term) | | 4 | |
| -1(3 term) | 1.3 | 5 | 2.74669E-13 |
| -1(4 term) | | 5 | 2.46692E-13 |
| -0.5(2 term) | | 4 | |
| -0.5(3 term) | 1.3 | 5 | 4.68514E-13 |
| -0.5(4 term) | | 5 | 5.18252E-13 |
| -1/3(2 term) | | 4 | |
| -1/3(3 term) | 1.3 | 5 | 5.50893E-13 |
| -1/3(4 term) | | 5 | 6.38156E-13 |
| 0.5(2 term) | | | 4.44089E-16 |
| 0.5(3 term) | 1.3 | 5 | 1.15286E-12 |
| 0.5(4 term) | | | 1.35025E-12 |
| 1.5(2 term) | | 6 | |
| 1.5(3 term) | 1.3 | 5 | 4.02306E-10 |
| 1.5(4 term) | | 5 | 1.54803E-11 |
| 2(2 term) | 1.3 | 6 | |
| 2(3 term)=N-method | 1.3 | 5 | 1.67788E-12 |

| Table 5.1.2 Caluculations of (2.5) in case of root 2 | | | |
|---|---|---|---|
| $q$ | $x_0$ | iterations | absolute error |
| N-method | 1.505 | 11 | 4.44089E-16 |
| (q=1) | 1.58 | 7 | 4.44089E-16 |
| | 2.27 | 4 | 1.74167E-11 |
| | 2.6 | 5 | |
| 0.5(2 term) | | 14 | |
| 0.5(3 term) | 1.505 | 16 | |
| 0.5(4 term) | | 23 | |
| 1.5(2 term) | | 6 | |
| 1.5(3 term) | 1.58 | 7 | 2.22045E-16 |
| 1.5(4 term) | | 7 | 2.22045E-16 |
| 2(2 term) | 2.27 | 4 | 2.37588E-13 |
| 2(3 term)=N-method | 2.27 | 4 | 1.74163E-11 |
| 3(2 term) | | 5 | |
| 3(3 term) | 1.58 | 6 | |
| 3(4 term)=N-method | | 7 | 4.44089E-16 |
| 4(2 term) | | 4 | |
| 4(3 term) | 1.58 | 6 | |
| 4(4 term) | | 6 | |
| 5(2 term) | | 4 | |
| 5(3 term) | 1.58 | 6 | |
| 5(4 term) | | 6 | |
| 6(2 term) | | 4 | |
| 6(3 term) | 1.58 | 6 | |
| 6(4 term) | | 6 | |
| 7(2 term) | | 5 | |
| 7(3 term) | 1.58 | 6 | |
| 7(4 term) | | 6 | |
| 8(2 term) | | 4 | |
| 8(3 term) | 1.58 | 6 | |
| 8(4 term) | | 5 | |
| 9(2 term) | | 4 | |
| 9(3 term) | 1.58 | 6 | |
| 9(4 term) | | 6 | |
| 10(2 term) | 2.27 | 6 | |
| 10(3 term) | 2.27 | 5 | |
| 10(4 term) | 2.6 | 6 | |

Table 5.2.1 Calculations of (5.2.1),(4.15),(4.17)

| $q$ | $|\mu_q(1^q)|$ | Right-hand side of (5.2.1) | $|\mu(1)|$ | $|1+f'(1)(1^-q)|/|1 \cdot f''(1)|$ | Right-hand side of (4.17) |
|---|---|---|---|---|---|
| -4 | 0.171201618 | 0.114134412 | 0.707106781 | 1.5 | 6.195386388 |
| -3.5 | 0.181419613 | 0.14513569 | 0.707106781 | 1.25 | 4.872039263 |
| -3 | 0.18973666 | 0.18973666 | 0.707106781 | 1 | 3.726779962 |
| -2.5 | 0.192098626 | 0.256131501 | 0.707106781 | 0.75 | 2.760717751 |
| -2 | 0.178885438 | 0.357770876 | 0.707106781 | 0.5 | 1.976423538 |
| -1.5 | 0.128007738 | 0.51203095 | 0.707106781 | 0.25 | 1.380984452 |
| -0.5 | 0.178885438 | 0.715541753 | 0.707106781 | 0.25 | 0.988211769 |
| 0.5 | 0.536656315 | 0.715541753 | 0.707106781 | 0.75 | 0.988211769 |
| 1 | 0.707106781 | 0.707106781 | 0.707106781 | 1 | 1 |
| 1.5 | 0.640038688 | 0.51203095 | 0.707106781 | 1.25 | 1.380984452 |

Table 5.2.2 Calculations of (5.2.2),(4.15),(4.17)

| $q$ | $|\mu_q(2^q)|$ | Right-hand side of (5.2.2) | $|\mu(2)|$ | $|1+f'(2)(1^-q)/2 \cdot f''(2)|$ | Right-hand side of (4.17) |
|---|---|---|---|---|---|
| -2 | 0.798940882 | 0.456537647 | 0.707106781 | 1.75 | 1.548846597 |
| -1 | 0.684806471 | 0.456537647 | 0.707106781 | 1.5 | 1.548846597 |
| 1 | 0.707106781 | 0.707106781 | 0.707106781 | 1 | 1 |
| 2 | 0.085600809 | 0.114134412 | 0.707106781 | 0.75 | 6.195386388 |
| 3 | 0.006872729 | 0.013745459 | 0.707106781 | 0.5 | 51.44293798 |
| 4 | 0.000487567 | 0.001950267 | 0.707106781 | 0.25 | 362.5691315 |
| 5 | 0 | 0.000312427 | 0.707106781 | 0 | 2263.272051 |
| 6 | 1.35628E-05 | 5.42513E-05 | 0.707106781 | 0.25 | 13033.92252 |
| 7 | 4.98242E-06 | 9.96485E-06 | 0.707106781 | 0.5 | 70960.11004 |
| 8 | 1.43051E-06 | 1.90735E-06 | 0.707106781 | 0.75 | 370728.1304 |
| 9 | 3.7676E-07 | 3.7676E-07 | 0.707106781 | 1 | 1876809.006 |
| 10 | 9.53674E-08 | 7.62939E-08 | 0.707106781 | 1.25 | 9268190.533 |
| 11 | 2.36448E-08 | 1.57632E-08 | 0.707106781 | 1.5 | 44858040.14 |

Shunji HORIGUCHI
315-0042
1-4-4 Baraki Ishioka Ibaragi, Japan
E-mail: shunhori@seagreen.ocn.ne.jp